\documentclass[11pt,a4paper]{article}
\usepackage{a4wide,amsfonts,amsmath,latexsym,amssymb,euscript,eufrak,graphicx,units,mathrsfs}

\usepackage{float}
\newfloat{figure}{H}{lof}
\floatname{figure}{\figurename}

\DeclareMathAlphabet{\eufrak}{U}{}{}{}  
\SetMathAlphabet\eufrak{normal}{U}{euf}{m}{n}
\SetMathAlphabet\eufrak{bold}{U}{euf}{b}{n}

\usepackage{amsmath, amsthm, amsfonts, amssymb}

\numberwithin{equation}{section}

\def\real{{\mathord{{\rm I\kern-2.8pt R}}}}        
\def\inte{{\mathord{{\rm I\kern-2.8pt N}}}}
\def\PP{{\mathord{{\rm I\kern-2.8pt P}}}}

\def\real{{\mathord{\mathbb R}}}

\def\inte{{\mathord{\mathbb N}}}

\newcommand{\retirer}[1]{$ $\newline  }

\def\HH{\EuFrak H}

\def\zit{\mathbb{Z}}

\def\E{\mathop{\hbox{\rm I\kern-0.20em E}}\nolimits}

\def\real{\mathbb{R}}

\newtheorem{prop}{Proposition}[section]

\newtheorem{theorem}[prop]{Theorem}
\newtheorem{thm}[prop]{Theorem}
\newtheorem{remark}[prop]{Remark}

\allowdisplaybreaks
\begin{document}

\begin{center}
{\Large{\bf Exact confidence intervals for the Hurst parameter of a fractional Brownian motion}}
\normalsize
\\~\\ by Jean-Christophe Breton\footnote{Universit\'e de La Rochelle, Laboratoire Math\'ematiques, Image et Applications, Avenue Michel Cr\'epeau, 17042 La Rochelle Cedex, France. Email: \texttt{jean-christophe.breton@univ-lr.fr}}, Ivan Nourdin\footnote{Laboratoire de Probabilit\'es et Mod\`eles Al\'eatoires, Universit\'e
 Pierre et Marie Curie (Paris VI), Bo\^ite courrier 188, 4 place Jussieu, 75252 Paris Cedex 05, France. Email: \texttt{ivan.nourdin@upmc.fr}},
and Giovanni Peccati\footnote{Equipe Modal'X, Universit\'{e} Paris Ouest -- Nanterre la D\'{e}fense, 200 Avenue de la R\'epublique, 92000 Nanterre, and LSTA, Universit\'{e} Paris VI, France. Email: \texttt{giovanni.peccati@gmail.com}}\\ {\it Universit\'e de La Rochelle, Universit\'e Paris VI and Universit\'e Paris Ouest}\\~\\
\end{center}

{\small \noindent {{\bf Abstract}: In this short note, we show how to use concentration inequalities in order to build exact confidence intervals for the Hurst parameter associated with a one-dimensional fractional Brownian motion.}
\\

\noindent {\bf Key words}: Concentration Inequalities; Exact confidence intervals; Fractional Brownian motion; Hurst parameter. \\

\noindent {\bf 2000 Mathematics Subject Classification:} 60F05; 60G15; 60H07.

\section{Introduction}

Let $B=\{B_t : t\geq 0\}$ be a fractional Brownian motion with Hurst index $H\in (0,1)$.
Recall that this means that $B$ is a real-valued continuous centered Gaussian process, with covariance given by $$E(B_tB_s) =\frac12(s^{2H}+t^{2H}-|t-s|^{2H}).$$
The reader is referred e.g. to \cite{PipTaqSurv} for a comprehensive introduction to fractional Brownian motion. We suppose that $H$ is unknown and verifies $H\leq H^*<1$, with $H^*$ known (throughout the paper, this is the only assumption we will make on $H$). Also, for a {\it fixed} $n\geq 1$, we assume that one observes $B$ at the times belonging to the set $\{k/n;\,k=0,\ldots,n+1\}$.

\smallskip

The aim of this note is to exploit the concentration inequality proved in \cite{NV}, in order to derive an \textsl{exact} (i.e., non-asymptotic) confidence interval for $H$. Our formulae hinge on the class of statistics
\begin{equation}\label{Intro1}
S_n=\sum_{k=0}^{n-1}\big(B_{\frac{k+2}n}-2B_{\frac{k+1}n}+B_{\frac kn}\big)^2, \quad n\geq 1.
\end{equation}
We recall that, as $n\rightarrow\infty$ and for every $H\in(0,1)$,
\begin{equation}\label{Intro2}
n^{2H-1}\,S_n \rightarrow 4-4^H, \,\, {\rm a.s. -} P,
\end{equation}
(see e.g. \cite{IL}), and also
\begin{eqnarray}\label{Intro3}
Z_n &=& n^{2H-\frac12}\,S_n -\sqrt{n}( 4-4^H) \\
&=& \frac{1}{\sqrt{n}}\sum_{k=0}^{n-1}\left(n^{2H}\big(B_{\frac{k+2}n}-2B_{\frac{k+1}n}+B_{\frac kn}\big)^2-(4-4^H)\right)\stackrel{\rm Law } {\Longrightarrow} N(0,c_H),\label{Intro5}
\end{eqnarray}
where $N(0,c_H)$ indicates a centered normal random variable, with finite variance $c_H>0$ depending only on $H$ (the exact expression of $c_H$ is not important for our discussion). We stress that the CLT (\ref{Intro5}) holds for every $H\in (0,1)$: this result should be contrasted with the asymptotic behavior of other remarkable statistics associated with the paths of $B$ (see e.g. \cite{BretNour} and \cite{BM}), whose asymptotic normality may indeed depend on $H$. The fact that $Z_n$ verifies a CLT for every $H$ is crucial in order to determine the asymptotic properties of our confidence intervals: see Remark \ref{REM fin} for further details.

\smallskip

The problem of estimating the self-similarity indices, associated with Gaussian and non-Gaussian stochastic processes, is crucial in applications, ranging from time-series, to physics and mathematical finance (see e.g. \cite{Nbook} for a survey). This issue has generated a vast literature: see \cite{Beran} and \cite{FoxTaqqu} for some classic references, as well as \cite{Coeurjolly}, \cite{GirRob}, \cite{IL}, \cite{TV}, and the references therein, for more recent discussions. However, the results obtained in our paper seems to be the first {\sl non-asymptotic} construction of a confidence interval for the Hurst parameter $H$. Observe that the knowledge of explicit non-asymptotic confidence intervals may be of great practical value, for instance in order to evaluate the accuracy of a given estimation of $H$ when only a fixed number of observations is available. In order to illustrate the novelty of our approach (i.e., replacing CLTs with concentration inequalities in the obte!
 ntion of confidence intervals), we also decided to keep things as simple as possible. In particular, we defer to a separate study the discussion of further technical points, such as e.g. the optimization of the constants appearing in our proofs.

\smallskip

The rest of this short note is organized as follows. In Section \ref{S : Mall} we state a concentration inequality that is useful for the discussion to follow. In Section \ref{S : Main} we state and prove our main result.

\section{A concentration inequality for quadratic forms}\label{S : Mall}

Consider a finite centered Gaussian family $X = \{X_k : k=0,...,M\}$, and write $R(k,l) = E(X_k X_l)$. In what follows, we shall consider two quadratic forms associated with $X$ and with some real coefficient $c$. The first is obtained by summing up the squares of the elements of $X$, and by subtracting the corresponding variances:
\begin{equation}\label{Quad1}
Q_1(c,X) = c\sum_{k=0}^M (X^2_k - R(k,k));
\end{equation}
the second quadratic form is
\begin{equation}\label{Quad2}
Q_2(c,X) = 2c^2\sum_{k,l=0}^M X_kX_l R(k,l).
\end{equation}
Note that $Q_2(c,X)\geq 0$. It is well known that, if $Q_1(c,X)$ is not a.s. zero, then the law of $Q_1(c,X)$ admits a density with respect to the Lebesgue measure (this claim can be easily proved by observing that $Q_1(c,X)$ can always be represented as a linear combination of independent centered $\chi^2$ random variables -- see \cite{Shik} for a general reference on similar results). The following statement, whose proof relies on the Malliavin calculus techniques developed in \cite{NV}, characterizes the tail behavior of $Q_1(c,X)$.
\begin{theorem}
\label{theo:NV}Let the above assumptions prevail, suppose that $Q_1(c,X)$ is not a.s. zero and fix $\alpha\geq 0$ and $\beta>0$.
Assume that $Q_2(c,X)\leq \alpha Q_1(c,X)+\beta$, a.s.-$P$. Then, for all $z>0$, we have
$$
P(Q_1(c,X)\geq z)\leq \exp\left(-\frac{z^2}{2\alpha z+2\beta}\right)\quad\mbox{and}\quad
P(Q_1(c,X)\leq -z)\leq \exp\left(-\frac{z^2}{2\beta}\right).
$$
In particular, $P(|Q_1(c,X)|\geq z)\leq 2\,\exp\left(-\frac{z^2}{2\alpha z+2\beta}\right).$
\end{theorem}
\noindent {\bf Proof.} In this proof, we freely use the language of isonormal Gaussian processes and Malliavin calculus; the reader is referred to \cite[Chapter 1]{Nbook} for any unexplained notion or result. Without loss of generality, we can assume that the Gaussian random variables $X_k$ have the form $X_k = X(h_k)$, where $X(\HH)=\{X(h) : h\in \HH\}$ is an isonormal Gaussian process over $\HH=\mathbb{R}^M$, and $\{h_k : k=1,...,M\}$ is a finite subset of $\HH$ verifying
$$
E[X(h_k)X(h_l)] = R(k,l) = \langle h_k,h_l\rangle_\HH.
$$
It follows that $Q_1(c,X) = I_2(c\sum_{k=0}^{M}h_k \otimes h_k)$, where $I_2$ stands for a double Wiener-It\^{o} stochastic integral with respect to $X$, so that the
$\HH$-valued Malliavin derivative of $Q_1(c,X)$ is given by
$$
DQ_1(c,X)=2c\sum_{k=0}^{M}X(h_k)h_k.
$$
Now write $L^{-1}$ for the pseudo-inverse of the Ornstein-Uhlenbeck generator associated with $X(\HH)$. Since $Q_1(c,X)$ is an element of the second Wiener chaos of $X(\HH)$, one has that $L^{-1}Q_1(c,X)=-\frac12\,Q_1(c,X)$. One therefore infers the relation
$$
\langle DQ_1(c,X),-DL^{-1}Q_1(c,X)\rangle_{\HH} = \frac12 \|DQ_1(c,X)\|^2_{\HH} = Q_2(c,X).
$$
The conclusion is now obtained by using the following general result.
\qed
\begin{thm}{\bf (See \cite[Theorem 4.1]{NV})}. Let $X(\HH)=\{X(h) : h\in \HH\}$ be an isonormal Gaussian process over some real separable Hilbert space $\HH$.
Write $D$ (resp. $L^{-1}$) to indicate the Malliavin derivative (resp. the pseudo-inverse of the generator $L$ of the Ornstein-Uhlenbeck semigroup).
Let $Z$ be a centered element of $\mathbb{D}^{1,2}:={\rm dom} D$, and suppose moreover that 
the law of $Z$ has a density with respect to the Lebesgue measure. 
If, for some $\alpha>0$ and $\beta\geq 0$, we have
$$
\langle DZ,-DL^{-1}Z\rangle_\HH\leq \alpha Z+\beta,\quad\mbox{a.s.-P},
$$
then, for all $z>0$, we have
$$
P(Z\geq z)\leq {\rm exp}\left(-\frac{z^2}{2\alpha z+2\beta}\right)\quad\mbox{and}\quad
P(Z\leq -z)\leq {\rm exp}\left(-\frac{z^2}{2\beta}\right).
$$
\end{thm}

\begin{remark}{\rm
One of the advantages of the concentration inequality stated in Theorem \ref{theo:NV} (with respect to other estimates that could be obtained by using the general inequalities by Borell \cite{borell}) is that they only involve {\sl explicit} constants.
}
\end{remark}

\section{Main result} \label{S : Main}

We go back to the assumptions and notation detailed in the Introduction. In particular, $B$ is a fractional Brownian motion with unknown Hurst parameter $H\in (0,H^*]$,
with $H^*<1$ known. The following result is the main finding of the present note.

\begin{thm}\label{intconf}
Fix $n\geq 1$, define $S_n$ as in (\ref{Intro1}) and fix a real $a$ such that
$0<a<(4-4^{H^*})\sqrt{n}$. For $x\in(0,1)$, set $g_n(x)=x-\frac{\log(4-4^x)}{2\log n}$. Then,
with probability at least
\begin{equation}\label{TOM}
\varphi(a) = \left[1-2\,{\rm exp}\left(
-\frac{a^2}{71\big(\frac{a}{\sqrt{n}}+3\big)}
\right)\right]_+,
\end{equation}
(where $[\cdot]_+$ stands for the positive part function), the unknown quantity $g_n(H)$ belongs to the following confidence interval:
$$
I(n) =[I_l(n), I_r(n)]= \left[
\frac12-\frac{\log S_n}{2\log n}+\frac{
\log\left(1-
\frac{a}{(4-4^{H^*})\sqrt{n}}
\right)
}{2\log n};
\frac12-\frac{\log S_n}{2\log n}+\frac{
\log\left(1+
\frac{a}{(4-4^{H^*})\sqrt{n}}
\right)
}{2\log n}
\right].
$$
\end{thm}

\smallskip

\begin{remark}{\rm
\begin{enumerate}
\item
We have that $\lim_{n\to\infty}g_n(H)=H$. Moreover, it is easily seen that the asymptotic relation (\ref{Intro2}) implies that, a.s.-$P$,
\begin{equation}\label{star}
\lim_{n\rightarrow \infty} I_l(n) = \lim_{n\rightarrow \infty} I_r(n) = H,
\end{equation}
that is, as $n \rightarrow \infty$, the confidence interval $I(n)$ ``collapses'' to the one-point set $\{H\}$.
\item In order to deduce (from Theorem \ref{intconf}) a genuine confidence interval for $H$, it is sufficient to (numerically) inverse the function $g_n$.
This is possible, since one has that $g_n'(x)\geq 1$ for every $x\in(0,1)$, thus yielding that $g_n$ is a continuous and strictly increasing bijection
from $(0,1)$ onto $(-\log 3/(2\log n),+\infty)$. It follows from Theorem \ref{intconf} that, with probability at least $\varphi(a)$, the parameter $H$ belongs to the interval
$$
J(n)=
[J_l(n),J_r(n)]=\left[g_n^{-1}\big(u(n) \big);
g_n^{-1}\big(I_r(n)\big)\right],
$$
where $u(n)= \max \{I_l(n); -\log 3/(2\log n)\}$. Observe that, since relation (\ref{star}) is verified, one has that $I_l(n)>-\log 3/(2\log n)$, a.s.-$P$, for $n$ sufficiently large.
Moreover, since $g_n^{-1}$ is $1$-Lipschitz, we infer that
\begin{eqnarray*}
J_r(n)-J_l(n)
&\leq& I_r(n)-I_l(n) = \frac{1}{2\log n}\,
\log\left(\frac
{(4-4^{H^*})\sqrt{n}+a}
{(4-4^{H^*})\sqrt{n}-a}
\right)
\end{eqnarray*}
so that, for every fixed $a$, the length of the confidence interval $J(n)$ converges a.s. to zero,
as $n\to\infty$, 
at the rate $O\big(1/(\sqrt{n}\log n)\big)$.
\item We now describe how to concretely build a confidence interval by means of Theorem \ref{intconf}. Start by fixing the error probability $\varepsilon$
(for instance, $\varepsilon=0,05$ or $0,01$). One has therefore two possible situations:

(i) If there are no restrictions on $n$ (that is, if the number of observations can be indefinitely increased), select first $a>0$
in such a way that
\begin{equation}\label{choice}
\exp\left(-\frac{a^2}{71(a+3)}\right)\leq \frac{\varepsilon}{2}
\end{equation}
(ensuring that $\varphi(a)\geq 1-\varepsilon$). Then, choose $n$ large enough in order to have
$$
\frac{a}{(4-4^{H^*})\sqrt{n}}< 1
\quad\mbox{and}\quad
\frac{1}{2\log n}\,
\log\left(\frac
{(4-4^{H^*})\sqrt{n}+a}
{(4-4^{H^*})\sqrt{n}-a}
\right)\leq L,
$$
where $L$ is some fixed (desired) upper bound for the length of the confidence interval.

(ii) If $n$ is fixed, then one has to select $a>0$
such that
$$
\exp\left(-\frac{a^2}{71\big(\frac{a}{\sqrt{n}}+3\big)}\right)
\leq \frac{\varepsilon}{2}\quad\mbox{and}\quad
a<(4-4^{H^*})\sqrt{n}.$$
If such an $a$ exists (that is, if $n$ is large enough), one obtains a confidence interval for $H$ 
of length less or equal to
$\frac{1}{2\log n}\,
\log\left(\frac
{(4-4^{H^*})\sqrt{n}+a}
{(4-4^{H^*})\sqrt{n}-a}
\right)$.

\item The fact that we work in a non-asymptotic framework is reflected by the necessity of choosing values of $a$ in such a way that the relation (\ref{choice}) is verified. On the other hand, if one uses directly the CLT (\ref{Intro5}) (thus replacing $Z_n$ with a suitable Gaussian random variable), then one can define an asymptotic confidence interval by selecting a value of $a$ such that a condition of the type $$\exp(-{\rm cst}\times a^2)\leq \varepsilon$$ is verified.
\item
By a careful inspection of the proof of Theorem \ref{intconf}, we see that the existence of $H^*$ is not required if we are only
interested in testing $H<\overline{H}$ for a given $\overline{H}$.
\end{enumerate}
}
\end{remark}

\smallskip

\noindent {\bf Proof of Theorem \ref{intconf}}. Define $X_n = \{X_{n,k} : k=0,...,n-1\}$, where
$$
X_{n,k} = B_{\frac{k+2}n}-2B_{\frac{k+1}n}+B_{\frac kn}.
$$
By setting
$$
\rho_H(r)=\frac12\big(
-|r-2|^{2H}+4|r-1|^{2H}-6|r|^{2H}+4|r+1|^{2H}-|r+2|^{2H}
\big), \quad r\in \mathbb{Z},
$$
one can prove by standard computations that the covariance structure of the Gaussian family $X_n$ is described by the relation $E(X_{n,k}X_{n,l}) = \rho_H (k-l)/n^{2H}$. Now let $Z_n$ be defined as in (\ref{Intro3}): it easily seen that $Z_n = Q_1(n^{2H-1/2},X_n)$ (as defined in (\ref{Quad1})).
We also have, see formula (\ref{Quad2}):
\begin{eqnarray}
Q_2(n^{2H-1/2},X_n)
&=&2n^{4H-1}\sum_{k,l=0}^{n-1} X_{n,k}X_{n,l}\frac{\rho_H(k-l)}{n^{2H}}\notag\\
&\leq&2n^{2H-1}\sum_{k,l=0}^{n-1} |X_{n,k}| |X_{n,l}| |\rho_H(k-l)|\notag\\
&\leq&n^{2H-1}\sum_{k,l=0}^{n-1} (X_{n,k}^2+X_{n,l}^2) |\rho_H(k-l)|\notag\\
&=&2n^{2H-1}\sum_{k,l=0}^{n-1} X_{n,k}^2 |\rho_H(k-l)|
\leq 2n^{2H-1}\sum_{k=0}^{n-1} X_{n,k}^2 \sum_{r\in\zit}|\rho_H(r)|\notag\\
&=&\frac{2}{\sqrt{n}}\left(\sum_{r\in\zit}|\rho_H(r)|\right)\left(Z_n+(4-4^H)\sqrt{n}\right)\notag\\
&\leq&\frac{2}{\sqrt{n}}\left(\sum_{r\in\zit}|\rho_H(r)|\right)\left(Z_n+3\sqrt{n}\right)
=\alpha_n Z_n+\beta\label{maj}
\end{eqnarray}
with
\begin{equation}\label{caravan}
\alpha_n=\frac2{\sqrt{n}}\sum_{r\in\zit}|\rho_H(r)|\quad\mbox{and}\quad
\beta=6\sum_{r\in\zit}|\rho_H(r)|.
\end{equation}
Since $Z_n \neq 0$, Theorem \ref{theo:NV} applies, yielding
\begin{equation}\label{bnd1}
P\big(|Z_n|>a\big)\leq 2\,{\rm exp}\left(
-\frac{a^2}{4\sum_{r\in\zit}|\rho_H(r)|\big(\frac{a}{\sqrt{n}}+3\big)}
\right).
\end{equation}
Now, let us find bounds on $\sum_{r\in\zit}|\rho_H(r)|$ that are independent of $H$.  Fix $r\geq 3$. Using
$$(1+u)^\alpha = 1+\sum_{k=1}^\infty \frac{\alpha(\alpha-1)\ldots(\alpha-k+1)}{k!}u^k\quad\mbox{for }\,0\leq u<1,$$
we can write
\begin{eqnarray*}
\rho_H(r)&=&\frac{r^{2H}}2\left(-\left(1-\frac2r\right)^{2H}+4\left(1-\frac1r\right)^{2H}-6+4\left(1+\frac1r\right)^{2H}-\left(1+\frac2r\right)^{2H}\right)\\
&=&\frac{r^{2H}}2\sum_{k=1}^{+\infty} \frac{2H(2H-1)\cdots(2H-k+1)}{k!}\big(-(-2)^k+4(-1)^k+4-2^k\big)r^{-k}\\
&=&r^{2H}\sum_{l=1}^{+\infty} \frac{2H(2H-1)\cdots(2H-2l+1)}{(2l)!}(4-4^l)r^{-2l}.
\end{eqnarray*}
Note that the sign of $2H(2H-1)\cdots(2H-2l+1)$ is the same as that of $2H-1$ and
\begin{eqnarray*}
\big|2H(2H-1)\cdots(2H-2l+1)\big|&=&2H\big|2H-1\big|(2-2H)\cdots (2l-1-2H)\\
&\leq& 2(2l-1)!.
\end{eqnarray*}
Hence,  we can write
\begin{eqnarray*}
|\rho_H(r)|
&\leq&r^{2H}\sum_{l=1}^{+\infty} \frac{4^l-4}{l}r^{-2l}\\
&=&4r^{2H}\log\left(1-\frac1{r^2}\right)
-r^{2H}\log\left(1-\frac4{r^2}\right)
\quad\big(\mbox{since $\log(1-u)=-\sum_{k=1}^{\infty} \frac{u^k}{k}$ if $0\leq u<1$}\big)\\
&\leq& \frac{243}{20}\,r^{2H-4}\quad\big(\mbox{since $4\log(1-u)-\log(1-4u)\leq\frac{243}{20}u^2$ if $0\leq u\leq\frac19$}\big)\\
&\leq& \frac{243}{20}\,r^{-2}.
\end{eqnarray*}
Consequently, taking into account of the fact that $\rho_H$ is an even function, we get
\begin{eqnarray}
\nonumber
\sum_{r\in\zit}|\rho_H(r)|
&\leq&|\rho_H(0)|+2|\rho_H(1)|+2|\rho_H(2)|+2\sum_{r= 3}^\infty|\rho_H(r)|\\
\nonumber
&=&|4-4^H|+|4\times 4^H-9^H-7|+|4-6\times 4^H+4\times 9^H-16^H|+2\sum_{r=3}^\infty|\rho_H(r)|\\
\nonumber
&\leq& 3+4+1+\frac{243}{10}\left(
\frac{\pi^2}6-1-\frac14
\right)= 17,59... \leq 17,75.
\label{eq:mu}
\end{eqnarray}
Putting this bound in (\ref{bnd1}) yields
\begin{equation}\label{bnd2}
P\big(|Z_n|>a\big)\leq 2\,{\rm exp}\left(
-\frac{a^2}{71\big(\frac{a}{\sqrt{n}}+3\big)}
\right).
\end{equation}
Note that the interest of this new bound is that the unknown parameter $H$ does not appear in the right-hand side.
Now we can construct the announced confidence interval for $g_n(H)$.
First, observe that $Z_n=n^{2H-\frac12}S_n- (4-4^H)\sqrt{n}$.
Using the assumption $H\leq H^*$ on the one hand, and (\ref{bnd2}) on the other hand, we get:
\begin{eqnarray*}
&&P\left(
\frac12-\frac{\log S_n}{2\log n}+\frac{
\log\left(1-
\frac{a}{(4-4^{H^*})\sqrt{n}}
\right)
}{2\log n}
\leq g_n(H)\leq
\frac12-\frac{\log S_n}{2\log n}+\frac{
\log\left(1+
\frac{a}{(4-4^{H^*})\sqrt{n}}
\right)
}{2\log n}
\right)\\
&\geq&P\left(
\frac12-\frac{\log S_n}{2\log n}+\frac{
\log\left(1-
\frac{a}{(4-4^{H})\sqrt{n}}
\right)
}{2\log n}
\leq H-\frac{\log(4-4^H)}{2\log n}\leq
\frac12-\frac{\log S_n}{2\log n}+\frac{
\log\left(1+
\frac{a}{(4-4^H)\sqrt{n}}
\right)
}{2\log n}
\right)\\
&=&
P\left(
\frac14-\frac{\log S_n}{2\log n}+\frac{
\log\left((4-4^{H})\sqrt{n}-a
\right)
}{2\log n}
\leq H\leq
\frac14-\frac{\log S_n}{2\log n}+\frac{
\log\left((4-4^{H})\sqrt{n}+a
\right)
}{2\log n}
\right)\\
&=&P\big(|Z_n|\leq a\big)\geq 1-2\,{\rm exp}\left(
-\frac{a^2}{71\big(\frac{a}{\sqrt{n}}+3\big)}
\right)
\end{eqnarray*}
which is the desired result. \qed

\medskip

\begin{remark} \label{REM fin}
{\rm The fact that $Q_2(n^{2H-1/2},X_n) \leq \alpha_n Z_n +\beta$ (see (\ref{maj})), where $\alpha_n \rightarrow 0$ and $\beta>0$, is consistent with the fact that $Z_n  \stackrel{\rm Law } {\Longrightarrow} N(0,c_H) $, and $Q_2(n^{2H-1/2},X_n) = \frac12\|DZ_n\|_\HH^2$, where $DZ_n$ is the Malliavin derivative of $Z_n$ (see the proof of Theorem \ref{theo:NV}). Indeed, according to Nualart and Ortiz-Latorre \cite{NO}, one has that $Z_n  \stackrel{\rm Law } {\Longrightarrow} N(0,c_H) $ if and only if $\frac12\|DZ_n\|_\HH^2$ converges to the constant $c_H$ in $L^2$. See also \cite{stein-ptrf} for a proof of this fact based on Stein's method.
}
\end{remark}

\medskip

\noindent {\bf Acknowledgment}. We are grateful to D. Marinucci for useful remarks.

\end{document}